\documentclass[12pt]{article}
\usepackage{enumerate}
\usepackage{amsmath,amssymb,amsthm, mathrsfs, mathtools}
\usepackage{latexsym}
\usepackage{graphics,graphicx}
\usepackage{hyperref}
\usepackage[bf, small]{titlesec}
\usepackage{authblk} %%% new line between authors
\usepackage{soul}
\usepackage{kotex}
\usepackage{lineno}
\theoremstyle{plain}
\newtheorem{Thm}{Theorem}[section]

\newtheorem{Lem}[Thm]{Lemma}
\newtheorem{Prop}[Thm]{Proposition}
\newtheorem{Cor}[Thm]{Corollary}

\theoremstyle{definition}

\newtheorem{Rmk}[Thm]{Remark}

\usepackage{standalone}
\usepackage{tikz}%pgfpages,,pgfkeys,pgfplots
\usetikzlibrary{arrows,positioning,matrix,fit,backgrounds,shapes,shapes.geometric, intersections}

\usetikzlibrary{shapes.callouts,decorations.pathmorphing} %%% 팝업용
\usepackage{calc} % calculate angles to evenly space vertices in circular arrangements.
\usetikzlibrary{decorations.markings} %%put arrowheads in the middle of directed edges
\tikzstyle{vertex}=[circle, draw, inner sep=0pt, minimum size=6pt] % style

\usetikzlibrary{arrows,matrix} %%% Hesse Diagram
\title{A digraph version of the Friendship Theorem}

\author[1]{Myungho Choi}
\author[1]{Hojin Chu}
\author[1]{Suh-Ryung Kim}

\affil[1]{\footnotesize Department of Mathematics Education, Seoul National University, Seoul 08826, Rep. of Korea}
\affil[ ]{\footnotesize\textit{nums8080@snu.ac.kr, ghwls8775@snu.ac.kr, srkim@snu.ac.kr}}
\date{}

\definecolor{LemonChiffon}{rgb}{100, 98, 80}
\definecolor{myblue}{rgb}{0,0.4,0.8}
\definecolor{orange}{rgb}{1, 0.4, 0}
\definecolor{mygreen}{rgb}{0, 0.8, 0}
\definecolor{myred}{rgb}{204, 0, 0}
\definecolor{violet}{RGB}{0.4,0.2,1}
\definecolor{brown}{rgb}{0.6, 0.4, 0}

 \newcounter{statement}

\begin{document}
\maketitle
\begin{abstract}
 The Friendship Theorem states that if in a party any pair of persons has precisely one common friend, then there is always a person who is everybody's friend and the theorem has been proved by Paul Erd\H{o}s, Alfr\'{e}d R\'{e}nyi, and Vera T. S\'{o}s in 1966. 
``What would happen if instead any pair of persons likes precisely one person?"
While a friendship relation is symmetric, a liking relation may not be symmetric.
Therefore to represent a liking relation we should use a directed graph.
We call this digraph a ``liking digraph".
It is easy to check that a symmetric liking digraph becomes a friendship graph if each directed cycle of length two is replaced with an edge.
In this paper, we provide a digraph formulation of the Friendship Theorem which characterizes the liking digraphs.
We also establish a sufficient and necessary condition for the existence of liking digraphs.
\end{abstract}



\section{Introduction}
In this paper, for graph-theoretical terminology and notations not defined, we follow \cite{bondy2010graph}.
In 1966, Erd\H{o}s {\it et al.} \cite{erd1966r} proved the Friendship Theorem that states if in a party any pair of persons has precisely one common friend, then there is always a person who is everybody's friend.
In graph terms, this theorem would be stated as follows: A finite graph in which each pair of vertices has exactly one common neighbor has a vertex adjacent to all the other vertices.
A graph is called a {\it friendship graph} if every pair of vertices has exactly one common neighbor.
This condition is said to be the {\it friendship condition}.

 We derive a digraph version of the Friendship Theorem.
 A digraph is called a {\it liking digraph} if every pair of its vertices has exactly one common out-neighbor.
The authors also introduce the notion of ``fancy wheel digraphs". 
A {\it fancy wheel digraph} is obtained from the disjoint union of directed cycles by adding one vertex $v$ with arcs to and from each vertex on the cycles. 
One may easily check that each fancy wheel digraph is a liking digraph
(see Figure~\ref{fig:fancy wheel graph} for an example).
\begin{figure}
\center	
\begin{tikzpicture}[scale=1.4]
     \tikzstyle{vertex} = [circle,draw,fill=white,inner sep=1.5];
      \draw[thick] (0,0) circle (0.5);
      \draw[thick] (0,0) circle (1);
      \draw[thick] (0,0) circle (1.5);
      \node[vertex] (0) at (0,0) [label=left:$v$] {};
      \node[vertex] (11) at (90:0.5) {};
      \draw[thick, stealth-stealth] (11) -- (0);
      \node[vertex] (12) at (-90:0.5) {};
      \draw[thick, stealth-stealth] (12) -- (0);
      \foreach \i in {1,2,3}{
      \node[vertex] (2\i) at (120*\i:1) {};
      \draw[thick, stealth-stealth] (2\i) -- (0);
      \path[tips,-stealth] (2\i) arc (120*\i:120*(\i+1)-3:1);
      }
      \foreach \i in {1,...,4}{
      \node[vertex] (3\i) at (45+90*\i:1.5) {};
      \draw[thick, stealth-stealth] (3\i) -- (0);
      \path[tips,-stealth] (3\i) arc (45+90*\i:45+90*(\i+1)-2:1.5);
      }
      \path[tips,-stealth] (11) arc (90:265:0.5);
      \path[tips,-stealth] (12) arc (-90:85:0.5);
    \end{tikzpicture}
 \caption{A fancy wheel digraph}
\label{fig:fancy wheel graph}
\end{figure}
A {\it $k$-diregular digraph} is a digraph in which each vertex has outdegree $k$ and indegree $k$ for a positive integer $k$.
A digraph is said to be {\it diregular} if it is a $k$-diregular digraph for some positive integer $k$.
Now, we are ready to state a digraph version of the Friendship Theorem.
\begin{Thm}\label{thm:friendship_digraph} (A digraph version of the Friendship Theorem)
If a graph is a liking digraph, then it is a fancy wheel digraph or it is $k$-diregular of order $k^2-k+1$ for some integer $k\ge 2$.
\end{Thm}
If we rephrase Theorem~\ref{thm:friendship_digraph} in terms corresponding to the original Friendship Theorem, it would be as follows:
	\begin{quote}
	If in a party any pair of persons likes precisely one person, then (i) there is a person who likes everybody and is liked by everybody or (ii) each person in the party likes exactly $k$ persons and is liked by exactly $k$ persons for some positive integer $k\geq 2$.
	Further, in the case (ii), there are $k^2-k+1$ people in the party.
\end{quote}

Since it is not common that a person likes oneself or to someone twice, the liking digraph in Theorem~\ref{thm:friendship_digraph} is assumed to be nontrivial, loopless, and without parallel arcs, and the same conditions are applied to the liking digraphs in this paper.

It is easy to check that a symmetric liking digraph becomes a friendship graph if each directed cycle of length two is replaced with an edge.
In this aspect, one can say that a liking digraph is a generalization of a friendship graph.
Actually, research that generalizes the concept of a friendship graph has been actively conducted.
For example,
Bose and Shrikhande~\cite{bose1969graphs} considered graphs in which each pair of vertices has
exactly $d$ common neighbors for a positive integer $d$. 
Delorme and Hahn~\cite{DELORME1984261} and Sudolsk\'{y}~\cite{sudolsky1978generalization} studied graphs in which any $t$ vertices have exactly $\lambda$ common neighbors for positive integers $t$ and $\lambda$.

We will also deduce the following interesting fact (Theorem~\ref{thm:reversed}).
 	\begin{quote}
			In a party, every pair of persons likes precisely one person if and only if every pair of persons is liked by precisely one person. 
	\end{quote}
	Based on the above fact, a liking digraph may be defined by a nontrivial loopless digraph without parallel arcs in which any two distinct vertices share exactly one in-neighbor. 
	
	By Theorem~\ref{thm:friendship_digraph},
	a liking digraph is a fancy wheel digraph or a diregular digraph.
	We note that the complete digraph of order $3$ is the only liking digraph that is both a $2$-diregular digraph and a fancy wheel digraph.
	While each fancy wheel digraph is a liking digraph, a $k$-diregular digraph of order $k^2-k+1$ may not be a liking digraph for some integer $k\ge 3$ (see Figure~\ref{fig:3-regular} for a $3$-diregular digraph of order $7$. It is not a liking digraph since $u$ and $v$ has no common out-neighbor). 
	Interestingly, the existence of a diregular liking digraph is determined by that of projective plane of certain order by the following theorem (the required details about projective plane are outlined in Section~\ref{sec:pre}).

\begin{Thm}\label{thm:SBIBD}
	For each integer $k\ge 3$, a $k$-diregular liking digraph exists if and only if there is a projective plane of order $k-1$.
\end{Thm}

The rest of this paper is devoted to proving the above two theorems.

\begin{figure}
\begin{center}
 \begin{tikzpicture}[auto,thick]
       \tikzstyle{player}=[minimum size=5pt,inner sep=0pt,outer sep=0pt,draw,circle]

    \tikzstyle{player1}=[minimum size=2pt,inner sep=0pt,outer sep=0pt,fill,color=black, circle]
    \tikzstyle{source}=[minimum size=5pt,inner sep=0pt,outer sep=0pt,ball color=black, circle]
    \tikzstyle{arc}=[minimum size=5pt,inner sep=1pt,outer sep=1pt, font=\footnotesize]
    \path (90:1.5cm)     node [player] [label=above:$u$] (v0){};
    \path (40:1.5cm)     node [player]  (v1){};
    \path (-10:1.5cm)     node [player]  (v2){};
    \path (-60:1.5cm)     node [player]  (v3){};
 \path (140:1.5cm)     node [player]  [label=above:$v$] (v6){};
  \path (190:1.5cm)     node [player]  (v5){};
    \path (240:1.5cm)     node [player]  (v4){};

   \draw[black,thick,-stealth] (v0) -> (v1);
 \draw[black,thick,-stealth] (v0) - +(v2);
  \draw[black,thick,-stealth] (v0) - +(v6);

  \draw[black,thick,-stealth] (v1) - +(v0);
  \draw[black,thick,-stealth] (v1) - +(v6);
  \draw[black,thick,-stealth] (v1) - +(v4);

  \draw[black,thick,-stealth] (v2) - +(v1);
  \draw[black,thick,-stealth] (v2) - +(v5);
  \draw[black,thick,-stealth] (v2) - +(v3);

  \draw[black,thick,-stealth] (v3) - +(v4);
  \draw[black,thick,-stealth] (v3) - +(v1);
  \draw[black,thick,-stealth] (v3) - +(v5);

  \draw[black,thick,-stealth] (v4) - +(v3);
  \draw[black,thick,-stealth] (v4) - +(v2);
  \draw[black,thick,-stealth] (v4) - +(v6);

  \draw[black,thick,-stealth] (v5) - +(v4);
  \draw[black,thick,-stealth] (v5) - +(v0);
  \draw[black,thick,-stealth] (v5) - +(v2);

  \draw[black,thick,-stealth] (v6) - +(v0);
  \draw[black,thick,-stealth] (v6) - +(v3);
  \draw[black,thick,-stealth] (v6) - +(v5);

    \end{tikzpicture}
 \end{center}
 \caption{A $3$-diregular digraph of order $7$ that is not a liking digraph}
\label{fig:3-regular}
\end{figure}
	
	\section{Preliminaries}\label{sec:pre}
	
We say that a digraph is {\it competitive} if any pair of its vertices has a common out-neighbor.
We note that every liking digraph is competitive.
The notion of competitive digraphs was introduced by Choi {\it et al.}~\cite{choi2022competitively}.
They showed that a tournament is competitive if and only if its order is greater than or equal to $7$.
They also completely characterized a competitive multipartite tournament in terms of the sizes of its partite sets. 
It is open to find other interesting families of competitive digraphs.

A {\it projective plane} consists of a set of objects called {\it points}, a second set of objects called {\it lines}, and a notion of when a point lies on a line, or equivalently, when a line passes through a point, so that any two lines have a unique point in common and any two points are incident with a unique line.

It is well-known that besides two degenerate examples, any finite projective plane satisfies the property that any two lines have the same number of points, that is, there is a number $q$ ($q\ge 2$) such that all lines of the plane have $q+1$ points. 
We call $q$ the {\it order} of the plane.
We note that the order of a projective plane is at least $2$.
Moreover, a projective plane of order $q$ satisfies that any point is incident with $q+1$ lines and both the number of points and the number of lines are $q^2+q+1$.
The two degenerate projective planes are the following.
\begin{enumerate}
\item[($P_1$)] There are a point $P$ and a line $l$ incident with $P$ such that all points are incident with $l$ and all lines are incident with $P$;
\item[($P_2$)] There are a point $P$ and a line $l$ not incident with $P$ such that every point except $P$ is incident with $l$ and every line except $l$ is incident with $P$.
\end{enumerate}
It is interesting to note that a fancy wheel digraph corresponds to a degenerate projective plane of type $(P_2)$.

%
%	Suppose that $A_1,A_2,\ldots,A_n$ are sets, which we refer to as a {\it set system}. A {\it (complete) system of distinct representatives} is a sequence $(a_1,a_2,\ldots,a_n)$ such that $a_i\in A_i$ for all $i$, and no two of the $a_i$ are the same.
%	Hall's Marriage Theorem plays a significant role in this paper, and it is as follows.
%\begin{Thm}[\cite{brualdi1977introductory}]\label{thm:MC}
%	The family $\mathcal{A}=(A_1,A_2,\ldots,A_n)$ of sets has a system of distinct representatives
%	if and only if 
%	\begin{center}
%	\begin{minipage}{0.85\textwidth}
%	(Hall's marriage condition)
%for each $1\leq k\leq n$ and each choice of $k$ distinct indices $i_1,i_2,\ldots,i_k$ from $\{1,2,\ldots,n\}$,
%\[\left|A_{i_1} \cup A_{i_2} \cup \cdots \cup A_{i_k}\right| \geq k.\]
%	\end{minipage}
%	\end{center}
%\end{Thm}

\section{Main results}
In this section,
we will focus on proving Theorems~\ref{thm:friendship_digraph} and \ref{thm:SBIBD}.
To this end, we show the following results first. 

\begin{Lem}\label{lem:degree at least 2}
	Each vertex of a liking digraph has outdegree at least $2$.
\end{Lem}
\begin{proof}
	Let $D$ be a liking digraph.
	Take a vertex $u$ of $D$.
	Since $D$ is not trivial, there is another vertex of $D$ and it has a common out-neighbor with $u$.
	Thus the outdegree of $u$ is at least one. 
	If a vertex $u$ has exactly one out-neighbor $v$, then $u$ and $v$ have no common out-neighbor since $D$ has no loops.
	Therefore the outdegree of $u$ is at least $2$.
	Since $u$ was arbitrarily chosen in $D$, the statement is true.
\end{proof}

\begin{Prop}\label{prop:equalities}
	Let $D$ be a liking digraph of order $n$.
	Then for each vertex $v \in V(D)$, the following equalities hold:
	\begin{enumerate}[(a)]
		\item $d^+(v)=d^-(v)$;
		\item $\sum_{u \in N^+(v)}(d^-(u)-1)=n-1$. 	\end{enumerate}
\end{Prop}
\begin{proof}
	Take a vertex $v$ in $D$.
By Lemma~\ref{lem:degree at least 2}, $v$ has outdegree at least $2$.
If $v$ has indegree at most $1$, then 
$d^-(v) \leq d^+(v)$.
Now suppose $v$ has indegree at least $2$ and take two distinct in-neighbors $x$ and $y$ of $v$.
Since any vertex distinct from $v$ has exactly one out-neighbor in $N^+(v)$ by the definition of liking digraphs, each of $x$ and $y$ has an out-neighbor in $N^+(v)$.
By the definition of liking digraphs again, $x$ and $y$ have distinct out-neighbors in $N^+(v)$ since $\{x,y\} \subseteq N^-(v)$.
Since $x$ and $y$ were arbitrarily chosen from $N^-(v)$, 
$d^-(v) \leq d^+(v)$.
Therefore we have shown that $d^-(v) \leq d^+(v) $ for each vertex $v$ in $D$. Thus
\[|A(D)|=\sum_{v\in V(D)}d^-(v) \leq \sum_{v\in V(D)}d^+(v)=|A(D)|\]
and so $d^-(v) =d^+(v)$ 
for each vertex $v$ in $D$.
Thus (a) holds.

We denote by $X$ a digraph isomorphic to the star on two edges with arcs pointing to the center.
Then by the definition of liking digraph, every pair of vertices are the leaves of exactly one copy of $X$.
Accordingly, a vertex is the leaf of exactly $n-1$ copies of $X$, that is, the number of copies of $X$ which have $v$ as a leaf equals to $n-1$ for every vertex $v$ in $D$.
Thus \[\sum_{u \in N^+(v)}(d^-(u)-1)=n-1. \qedhere\]
\end{proof}	

By Lemma~\ref{lem:degree at least 2} and Proposition~\ref{prop:equalities}, the following result is immediately true.

\begin{Cor}\label{cor:wheel}
	If there is a vertex $v$ of a liking digraph $D$ such that $N^+(v)=V(D)-\{v\}$, then each of the other vertices must have outdegree $2$.
\end{Cor}

Given a digraph $D$, $D^\leftarrow$ is the digraph obtained from $D$ by reversing the direction of each arc.
That is, $V(D^\leftarrow)=V(D)$ and \[(u,v) \in A(D^\leftarrow) \iff (v,u) \in A(D).\]

\begin{Thm}\label{thm:reversed}
	Let $D$ be a liking digraph.
	Then $D^{\leftarrow}$ is a liking digraph.
\end{Thm}
\begin{proof}
	If there are two common in-neighbors of $u$ and $v$ for some distinct vertices $u,v \in V(D)$, then the two common in-neighbors have $u$ and $v$ as two common out-neighbors, which violates the definition of liking digraph.
	Thus for each pair of vertices $u$ and $v$ in $V(D)$,
	\[|N^-(u)\cap N^-(v)| \le 1. \]
	One may check that
	\begin{align*}
		 \sum_{\{u,v\}\subseteq V(D)}1 &\ge 
		 \sum_{\{u,v\}\subseteq V(D)} |N^-(u)\cap N^-(v)| \\
		 &=\sum_{w\in V(D)}\sum_{\{u,v\}\subseteq N^+(w)}1\\
		&=\sum_{w\in V(D)}\sum_{\{u,v\}\subseteq N^-(w)}1\\
		&=\sum_{\{u,v\}\subseteq V(D)} |N^+(u)\cap N^+(v)| \\
		&=\sum_{\{u,v\}\subseteq V(D)} 1
	\end{align*}
	where the second equality is true by Proposition~\ref{prop:equalities}(a).
	Therefore, for each pair of vertices $u$ and $v$ in $V(D)$, 
	\[|N^-(u)\cap N^-(v)| = 1. \]
	Thus $D^{\leftarrow}$ is a liking digraph.
\end{proof}
 
Given a digraph $D$,
the maximum outdegree of $D$ is denoted by $\Delta^+(D)$.

\begin{Prop}\label{prop:outdegree_n-1}
	Let $D$ be a liking digraph of order $n$.
	Then $\Delta^+(D)= n-1$ if and only if $D$ is a fancy wheel digraph.
\end{Prop}
\begin{proof}
The ``if" part is immediately true by the definition of fancy wheel digraph.
To show the ``only if" part, suppose $\Delta^+(D)= n-1$.
That is, there is a vertex $v$ with $d^+(v)=n-1$.
Then $N^+(v)=V(D)-\{v\}$.
Thus, by Corollary~\ref{cor:wheel}, each vertex in $N^+(v)$ has indegree exactly $2$. 
Since $N^+(v)=V(D)-\{v\}$, each vertex of $D$ except $v$ has exactly one in-neighbor distinct from $v$.
Thus $D-v$ is a vertex-disjoint union of directed cycles.
Hence $D$ is a fancy wheel digraph.
\end{proof}

\begin{Lem}\label{lem:nice_inequality}
	Let $D$ be a liking digraph of order $n$.
	Then $(d^-(u)-1)(d^-(v)-1)\leq n-2$
	for any pair of vertices $u$ and $v$ in $D$.
\end{Lem}

\begin{proof}
Take two distinct vertices $u$ and $v$.
Then $d^+(u)=d^-(u)$ and $d^+(v)=d^-(v)$ by Proposition~\ref{prop:equalities}(a).
By Theorem~\ref{thm:reversed},
$u$ and $v$ have exactly one common in-neighbor $w$.
Let $N^-(u)=\{u_1,\ldots,u_{d^-(u)-1},w\}$ and $N^-(v)=\{v_1,\ldots,v_{d^-(v)-1},w\}$.
By the definition of liking digraph, for any integers $1\leq i \leq d^-(u)-1$ and $1\leq j \leq d^-(v)-1$,
$u_i$ and $v_j$ have exactly one common out-neighbor $w_{ij}$ in $D$.
Since $N^-(u) \cap N^-(v)=\{w\}$, $w_{ij} \neq u$ and $w_{ij} \neq v$ for any integers $1\leq i \leq d^-(u)-1$ and $1\leq j \leq d^-(v)-1$.
Thus \begin{equation}\label{eq:prop:nice_inequality}
\{w_{ij} \colon\, 1\leq i \leq d^-(u)-1,\ 1\leq j \leq d^-(v)-1\}\subseteq V(D)-\{u,v\}.
\end{equation}
If $w_{ij}=w_{i'j'}$ for some distinct integers $i$ and $i'$, then $u$ and $w_{ij}$ are common out-neighbors of $u_i$ and $u_{i'}$, which is impossible.
If $w_{ij}=w_{i'j'}$ for some distinct integers $j$ and $j'$, then $v$ and $w_{ij}$ are common out-neighbors of $v_j$ and $v_{j'}$, which is impossible.
Thus \[\left|\{w_{ij} \colon\, 1\leq i \leq d^-(u)-1,\ 1\leq j \leq d^-(v)-1\}\right|=(d^-(u)-1)(d^-(v)-1) \]
and so, by equation \eqref{eq:prop:nice_inequality},
$(d^-(u)-1)(d^-(v)-1)\leq n-2.$
\end{proof}

The following proposition plays a key role in proving Theorem~\ref{thm:friendship_digraph}.
\begin{Prop}\label{prop:nice_lemma}
	Let $D$ be a liking digraph.
	If $(u,v) \notin A(D)$ or $(v,u) \notin A(D)$ for some two distinct vertices $u$ and $v$,
	then $d^+(u)=d^+(v)$.
	\end{Prop}
	\begin{proof}
	Suppose that $u$ and $v$ are two vertices such that
	$(u,v) \notin A(D)$ or $(v,u)\notin A(D)$.
Without loss of generality,
we may assume \[(v,u) \notin A(D).\]
Let $n$ be the order of $D$.
Then 
\begin{align*}
(d^+_{D}(u)-1)(n-1)&=(d^-_{D}(u)-1)(n-1) \tag{by Proposition~\ref{prop:equalities}(a)}
\\&=\sum_{x\in N^+(v)}(d^-_{D}(u)-1)(d^-_{D}(x)-1) \tag{by Proposition~\ref{prop:equalities}(b)} 
\\& \le \sum_{x\in N^+(v)}(n-2) \tag{by Lemma~\ref{lem:nice_inequality} since $u\notin N^+(v)$} 
\\& = d^+_{D}(v)(n-2).  
\end{align*}
Thus $(d^+_{D}(u)-1)(n-1)< d^+_{D}(v)(n-1)$ and so $d^+_{D}(u)-1<d^+_{D}(v)$.
Therefore $d^+_{D}(u)\le d^+_{D}(v)$.
Now we consider the digraph $D^\leftarrow$.
Then $D^{\leftarrow}$ is a liking digraph of order $n$ by Theorem~\ref{thm:reversed} and $(u,v) \notin A(D^{\leftarrow})$.
Thus, by applying the above argument to $D^{\leftarrow}$,
we obtain $d^+_{D^{\leftarrow}}(v)\leq d^+_{D^{\leftarrow}}(u)$.
Then, by Proposition~\ref{prop:equalities}(a),
$d^-_{D^{\leftarrow}}(v)\leq d^-_{D^{\leftarrow}}(u)$ and so
$d^+_{D}(v)\leq d^+_{D}(u)$.
Therefore $d^+_{D}(u)=d^+_{D}(v)$.
\end{proof}

Now we are ready to prove Theorem~\ref{thm:friendship_digraph}.

\begin{proof}[Proof of Theorem~\ref{thm:friendship_digraph}]
	Let $D$ be a liking digraph of order $n$.
	If $\Delta^+(D)=n-1$, then $D$ is a fancy wheel digraph by Proposition~\ref{prop:outdegree_n-1}.
	Suppose $\Delta^+(D)< n-1$.
	Take a vertex $u$ and let $d^+(u)=k$ for some positive integer $k$.
	Then $k\ge 2$ by Lemma~\ref{lem:degree at least 2}.
	Since $\Delta^+(D)< n-1$, there is a vertex $v$ such that $(u,v)\notin A(D)$.
	By Proposition~\ref{prop:nice_lemma},
	\[ d^+(v)=d^+(u)=k.\]
	By the definition of liking digraph,
	$u$ and $v$ have exactly one common out-neighbor $w$.
	By Lemma~\ref{lem:degree at least 2}, $\Delta^+(D)\ge 2$ and so $n >3$.
	Then we may take a vertex $x$ in $V(D)-\{u,v,w\}$.
	Since $N^+(u)\cap N^+(v)=\{w\}$, $x\not\in N^+(u)$ or $x\not\in N^+(v)$.
	Thus $d^+(x)=d^+(u)$ or $d^+(x)=d^+(v)$ by Proposition~\ref{prop:nice_lemma}.
	Therefore \[d^+(x)=k.\]
	Since $x$ was arbitrarily chosen from $V(D)-\{u,v,w\}$,
	each vertex except $w$ has outdegree $k$.
	Since $\Delta^+(D)<n-1$,
	there exists a vertex $y$ such that $y \not\in N^+(w)$.
	Thus $d^+(y)=d^+(w)$ by Proposition~\ref{prop:nice_lemma}.
	Since $y \in V(D)-\{w\}$, $d^+(y)=k$ and so $d^+(w)=k$.
Hence each vertex in $D$ has outdegree $k$.
Therefore $D$ is a $k$-diregular digraph by Proposition~\ref{prop:equalities}(a).
Moreover, $n=k^2-k+1$ by Proposition~\ref{prop:equalities}(b).
\end{proof}

Now we prove Theorem~\ref{thm:SBIBD}.

\begin{proof}[Proof of Theorem~\ref{thm:SBIBD}]
		To show the ``only if" part, suppose that there is a $k$-diregular liking digraph $D$ for some integer $k \ge 3$.
	Then $D$ has order $k^2-k+1$ by Theorem~\ref{thm:friendship_digraph}.
	We consider the set $\mathcal{P}$ of points each of which is a vertex in $D$ and the set $\mathcal{L}$ of lines each of which is the in-neighborhood of a vertex.
	Since $D$ is $k$-diregular, each point in $\mathcal{P}$ is incident with exactly $k$ lines and each line in $\mathcal{L}$ is incident with exactly $k$ points.
	Moreover, since $D$ is a liking digraph, any two points in $\mathcal{P}$ are incident with a unique line.
	Indeed, since $D^{\leftarrow}$ is a liking digraph by Theorem~\ref{thm:reversed}, any two lines in $\mathcal{L}$ have a unique point in common.
	Therefore $\mathcal{P}$ and $\mathcal{L}$ form a projective plane of order $k-1$ and this proves the ``only if" part.

	To show the ``if" part, suppose that there is a projective plane $\mathbf{P}$ of order $k-1$ for some positive integer $k\ge 3$.
	Let $\mathcal{P}$ be the set of its points and $\mathcal{L}$ be the set of its lines. 
	Then  $|\mathcal{P}|=|\mathcal{L}|=k^2-k+1=:n$.
	We let $\mathcal{L}=\{l_1, \ldots, l_{n}\}$.
	Let $\mathcal{P}_i$ denote the points on the line $l_i$ for each $1\leq i\leq n$.
	Then $|\mathcal{P}_i|=k$ for each  $1\leq i\leq n$.
	It is easy to check that the collection $\{\mathcal{P}\setminus \mathcal{P}_1, \ldots, \mathcal{P}\setminus \mathcal{P}_n\}$ has a system of distinct representatives. 
	Let $\{P_1, \ldots, P_n\}$ be a system of distinct representatives of $\{\mathcal{P}\setminus \mathcal{P}_1, \ldots, \mathcal{P}\setminus \mathcal{P}_n\}$
	Now, we consider the digraph $D$ with the vertex set $\mathcal{P}$ and the arc set
	\[\bigcup_{t=1}^{n} \{(P,P_t) \colon\, P\in \mathcal{P}_t\}.\]
	Since $P_t\not\in \mathcal{P}_t$ for each $1\le t \le n$, $D$ is loopless.
	Further, since any two points are incident with a unique line, each pair of vertices in $D$ has a unique common out-neighbor.
	Thus $D$ is a liking digraph of order $k^2-k+1$. 
	Since each point is incident with exactly $k$ lines and each line has exactly $k$ points, $D$ is a $k$-diregular digraph.
	\end{proof}

\begin{Rmk}
	When we prove Theorem~\ref{thm:SBIBD}, we correspond a projective plane to a $k$-diregular liking digraph $D$.
	It is easy to check that, by construction, $D$ is a fancy wheel digraph if and only if the corresponding projective plane is degenerate and of type $(P_2)$. 
\end{Rmk}

\section{Concluding Remarks}
	While the classification of all projective planes is far from being complete, it is known that a projective plane of order $n$ exists if $n$ is a prime power.
	Moreover, there are results that a projective plane of order $n$ does not exist for $n\in \{6,10\}$ (see \cite{bruck1949planeoforder6, Lam1989planeoforder10}).
	Therefore the following result is obtained by Theorem~\ref{thm:SBIBD}:
	\begin{itemize}
		\item A liking digraph of order $k^2-k+1$ exists where $k=p^m+1$ for any prime $p$ and any positive integer $m$.
		\item A liking digraph of order $k^2-k+1$ does not exist for $k\in \{7,11\}$.
	\end{itemize}	

\section{Acknowledgement}
We are thankful to the reviewer whose suggestions highly advanced the clarity of the paper.

This work was supported by Science Research Center Program through the National Research Foundation of Korea(NRF) grant funded by the Korean Government (MSIT)(NRF-2022R1A2C1009648 and 2016R1A5A1008055).
%\bibliographystyle{abbrv}
%\bibliography{competition}

\end{document}